%% file: loclhypnothyp.tex
\documentclass[12pt,notitlepage, leqno]{amsart}

\include{packages_example}

\usepackage{subfig}
\usepackage{floatrow}

\makeatletter 
\def\paragraph{\@startsection{paragraph}{4}%
\z@\z@{-\fontdimen2\font}%
{\normalfont\bfseries}}
\makeatother

\begin{document}

	 	\title[]{A locally hyperbolic 3-manifold that is not hyperbolic}
	 	\author{Tommaso Cremaschi}
	 	\date{\today}
		\thanks{The author gratefully acknowledges support from the U.S. National Science Foundation grants DMS 1107452, 1107263, 1107367 "RNMS: GEometric structures And Representation varieties" (the GEAR Network) and also from the grant DMS-1564410: Geometric Structures on Higher Teichm\"uller Spaces.}

	 	\maketitle
	 	
		\small
	 	
		\paragraph{Abstract:} We construct a locally hyperbolic 3-manifold $M_\infty$ such that $\pi_ 1(M_\infty)$ has no divisible subgroup. We then show that $M_\infty$ is not homeomorphic to any complete hyperbolic manifold. This answers a question of Agol \cite{DHM,Ma2007}.
		\normalsize

	 	\begin{center}
	 		\section*{Introduction}
	 	\end{center}
		Throughout this paper, $M$ is always an oriented, aspherical 3-manifold. A 3-manifold $M$ is \textit{hyperbolizable} if its interior is homeomorphic to $\quotient{\mathbb H^3}\Gamma$ for $\Gamma\subgroup \text{Isom}(\mathbb H^3)$ a discrete, torsion free subgroup. An irreducible 3-manifold $M$ is of \textit{finite-type} if $\pi_1(M)$ is finitely generated and we say it is of \textit{infinite-type} otherwise.\ By Geometrization (2003, \cite{Per2003.1,Per2003.2,Per2003.3}) and Tameness (2004, \cite{AG2004,CG2006}) a finite type 3-manifold $M$ is hyperbolizable if and only if $M$ is the interior of a compact 3-manifold $\overline M$ that is atoroidal and with non finite $\pi_1(\overline M)$. On the other hand, if $M$ is of infinite type not much is known and we are very far from a complete topological characterisation.\ Nevertheless, some interesting examples of these manifolds have been constructed in \cite{SS2013,BMNS2016}. What we do know are necessary condition for a manifold of infinite type to be hyperbolizable. If $M$ is hyperbolizable then $M\cong \hyp\Gamma$, hence by discreteness of $\Gamma$ and the classification of isometries of $\bH$  we have that no element $\gamma\in\Gamma$ is \textit{divisible} (\cite[Lemma 3.2]{Fr2011}). Here, $\gamma\in\Gamma$ is divisible if there are infinitely many $\alpha\in\pi_1(M)$ and $n\in\N$ such that: $\gamma= \alpha^n$. We say that a manifold $M$ is \textit{locally hyperbolic} if every cover $N\twoheadrightarrow M$ with $\pi_1(N)$ finitely generated is hyperbolizable. Thus, local hyperbolicity and having no divisible subgroups in $\pi_1$ are necessary conditions. In \cite{DHM,Ma2007} Agol asks whether these conditions could be sufficient for hyperbolization:
		 
		 		\vspace{0.3cm}
		
			 	\begin{con2}[Agol] Is there a 3-dimensional manifold $M$ with no divisible elements in $\pi_1(M)$ that is locally hyperbolic but not hyperbolic?
	 	\end{con2}
We give a positive answer:
		
				\begin{customthm}{1}\label{agolquest}
		There exists a locally hyperbolic 3-manifold with no divisible subgroups in its fundamental group that does not admit any complete hyperbolic metric.
		\end{customthm}	
		
		\bigskip
 
		\paragraph*{Outline of the proof:}		
		The manifold $M_\infty$ is a thickening of the 2-complex obtained by gluing to an infinite annulus $A$ countably many copies of a genus two surface $\set{\Sigma_i}_{i\in\Z}$ along a fixed separating curve $\gamma$ such that the $i$-th copy $\Sigma_i$ is glued to $S^1\times \set i$. The manifold $M_\infty$ covers a compact non-atoroidal manifold $M$ containing an incompressible two sided surface $\Sigma$. Since $\pi_1(M_\infty)\subgroup \pi_1(M)$ and $M$ is Haken by \cite{Sh1975} we have that $\pi_1(M_\infty)$ has no divisible elements. By construction $M_\infty$ has countably many embedded genus two surfaces $\set{\Sigma_i}_{i\in\Z}$ that project down to $\Sigma$. By a surgery argument it can be shown that $M_\infty$ is atoroidal. Moreover, if we consider the lifts $\Sigma_{-i},\Sigma_i$ they co-bound a submanifold $M_i$ that is hyperbolizable and we will use the $M_i$ to show that $M_\infty$ is locally hyperbolic (see Lemma \ref{locallyhyperbolic}). Thus, $M_\infty$ satisfies the conditions of Agol's question.

		The obstruction to hyperbolicity arises from the lift $A$ of the essential torus $T$. The lift $A$ is an open annulus such that the intersection with all $M_i$ is an embedded essential annulus $A_i\eqdef A\cap M_i$ with boundaries in $\Sigma_{\pm i}$. The surfaces $\Sigma_{\pm i}$ in the boundaries of the $M_i$ have the important property that they have no homotopic essential subsurfaces except for the one induced by $A$. This gives us the property that both ends of $A$ see an 'infinite' amount of topology. This is in sharp contrast with finite type hyperbolic manifolds in which, by Tameness, every such annulus only sees a finite amount of topology.

		\medskip
				
		In future work we will give a complete topological characterisation of hyperbolizable 3-manifolds for a class of infinite type 3-manifolds. This class contains $M_\infty$ and the example of Souto-Stover \cite{SS2013} of a hyperbolizable Cantor set in $S^3$.

				\medskip
	
	 	\paragraph*{Acknowledgements:}: {I would like to thank J.Souto for introducing me to the problem and for his advice without which none of this work would have been possible. I would also like to thank I.Biringer and M.Bridgeman for many helpful discussions and for looking at some \textit{unreadable} early drafts. I am also grateful to the University of Rennes 1 for its hospitality while this work was done. }
	
\medskip

		\paragraph*{Notation:} We use $\simeq$ for homotopic and by $\pi_0(X)$ we intend the connected components of $X$. With $\Sigma_{g,k}$ we denote the genus $g$ orientable surface with $k$ boundary components. By $N\hookrightarrow M$ we denote embeddings while $S\imm M$ denotes immersions.
		
		\vspace{1cm}

		\section{Background} We now recall some fact and definitions about the topology of 3-manifolds, more details can be found in \cite{He1976,Ha2007,Ja1980}.

An orientable 3-manifold $M$ is said to be \emph{irreducible} if every embedded sphere $S^2$ bounds a 3-ball. A map between manifolds is said to be \emph{proper} if it sends boundaries to boundaries and pre-images of compact sets are compact. We say that a connected properly immersed surface $S\imm M$ is \emph{$\pi_1$-injective} if the induced map on the fundamental groups is injective. Furthermore, if $S\hookrightarrow M$ is embedded and $\pi_1$-injective we say that it is \emph{incompressible}. If $S\hookrightarrow M$ is a non $\pi_1$-injective two-sided surface by the Loop Theorem we have that there is a compressing disk $D\hookrightarrow M$ such that $\partial D=D\cap S$ and $\partial D$ is non-trivial in $\pi_1(S)$. 

An irreducible 3-manifold $(M,\partial M)$ is said to have \emph{incompressible} \emph{boundary} if every map: $(D^2,\partial D^2)\hookrightarrow (M,\partial M)$ is homotopic via a map of pairs into $\partial M$. Therefore, $(M,\partial M)$ has incompressible boundary if and only if each component $S\in\pi_0(S)$ is incompressible, that is $\pi_1$-injective.  An orientable, irreducible and compact $3$-manifold is called \textit{Haken} if it contains a two-sided $\pi_1$-injective surface. A 3-manifold is said to be \emph{acylindrical} if every map $(S^1\times I,\partial (S^1\times I))\rar (M,\partial M)$ can be homotoped into the boundary via maps of pairs.

\bdefi A 3-manifold $M$ is said to be \emph{tame} if it is homeomorphic to the interior of a compact 3-manifold $\overline M$. \edefi

Even 3-manifolds that are homotopy equivalent to compact manifolds need not to be tame.\ For example the Whitehead manifold \cite{Wh1935} is homotopy equivalent to $\R^3$ but is not homeomorphic to it.

\bdefi We say that a codimension zero submanifold $N\overset{\iota}{\hookrightarrow} M$ forms a \emph{Scott core} if the inclusion map $\iota_*$ is a homotopy equivalence. \edefi

		By \cite{Sc1973,HS1996,RS1990} given an orientable irreducible 3-manifold $M$ with finitely generated fundamental group a Scott core exists and is unique up to homeomorphism.		
		
Let $M$ be a tame 3-manifold, then given a Scott core $C\hookrightarrow M\subset \overline M$ with incompressible boundary we have that, by Waldhausen's cobordism Theorem \cite{Wa1968}, every component of $\overline{\overline M\setminus C}$ is a product submanifold homeomorphic to $S\times I$ for $S\in\pi_0(\partial C)$.

		\bdefi Given a core $C\hookrightarrow M$ we say that an end $E\subset\overline{ M\setminus C}$ is \textit{tame} if it is homeomorphic to $S\times [0,\infty)$ for $S=\partial E$.\edefi
		
		 A core $C\subset M$ gives us a bijective correspondence between the ends of $M$ and the components of $\partial C$. We say that a surface $S\in\pi_0(\partial C)$ \textit{faces the end $E$} if $E$ is the component of $\overline{M\setminus C}$ with boundary $S$. It is a simple observation that if an end $E$ facing $S$ is exhausted by submanifolds homeomorphic to $S\times I$ then it is a tame end.

		\section{Proof of Theorem 1}

Consider a surface of genus two $\Sigma$ and denote by $\alpha$ a separating curve that splits it into two punctured tori. To $\Sigma\times I$ we glue a thickened annulus $C\eqdef (S^1\times I)\times I$ so that $S^1\times I\times \set i$ is glued to a regular neighbourhood of $\alpha\times i$, for $i=0,1$. We call the resulting manifold $M$:
\begin{center}\begin{figure}[h!]
					
						\def\svgwidth{200pt}
						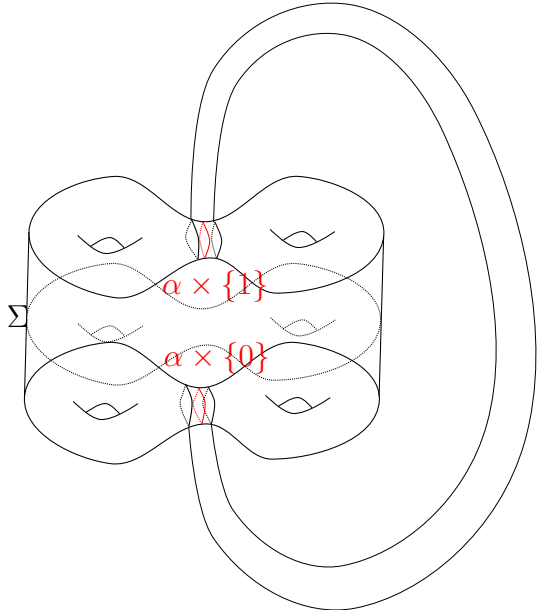
						\caption{The manifold $M$.}						
						\end{figure}
	 								\end{center}

The manifold $M$ is not hyperbolic since it contains an essential torus $T$ coming from the cylinder $C$. Moreover, $M$ has a surjection $p$ onto $S^1$ obtained by projecting the surfaces in $\Sigma\times I$ onto $I$ and also mapping the cylinder onto an interval. We denote by $H$ the kernel of the surjection map $p_*:\pi_1(M)\twoheadrightarrow\pi_1(S^1)$. 

Consider an infinite cyclic cover $M_\infty$ of $M$ corresponding to the subgroup $H$. The manifold $M_\infty$ is an infinite collection of $\set{\Sigma\times I}_{i\in\Z}$ glued to each other via annuli along the separating curves $\alpha\times\set {0,1}$. Therefore, we have the following covering:

		\begin{center}\begin{figure}[h!]
					
						\def\svgwidth{400pt}
						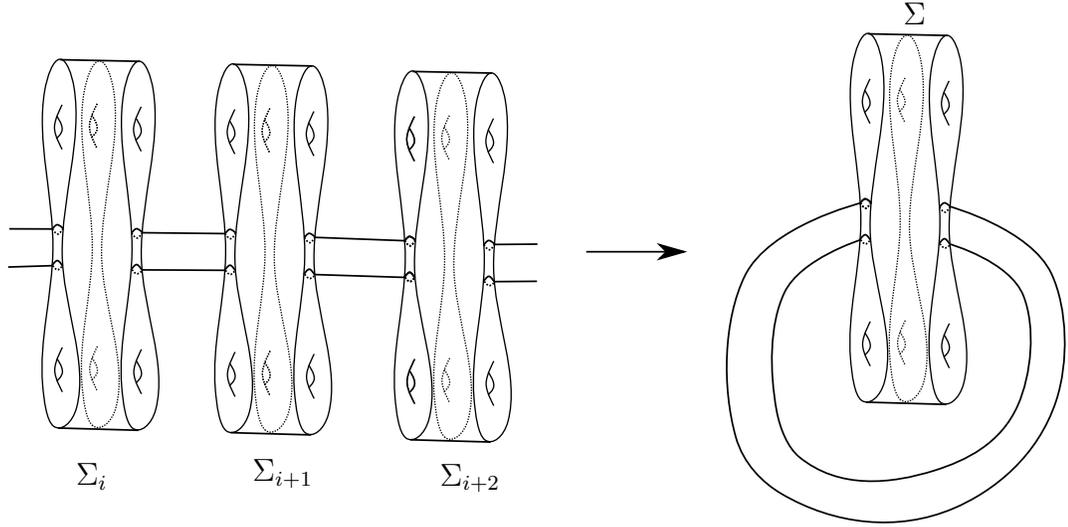
						\caption{The infinite cyclic cover.}
						\end{figure}
	 								\end{center}
					where the $\Sigma_i$ are distinct lifts of $\Sigma$ and so are incompressible in $M_\infty$. Since $\pi_1(M_\infty)$ is a subgroup of $\pi_1(M)$ and $M$ is Haken ($M$ contains the incompressible surface $\Sigma$) by \cite{Sh1975} we have that $\pi_1(M)$ has no divisible elements, thus $\pi_1(M_\infty)$ has no divisible subgroups as well.

					\blem\label{locallyhyperbolic} The manifold $M_\infty$ is locally hyperbolic. \elem\bpf

We claim that $M_\infty$ is atoroidal and exhausted by hyperbolizable manifolds. Let $T^2\hookrightarrow M_\infty$ be an essential torus with image $T$. Between the surfaces $\Sigma_i$ and $\Sigma_{i+1}$ we have incompressible annuli $C_i$ that separate them. Since $T$ is compact it intersects at most finitely many $\set{C_i}$. Moreover, up to isotopy we can assume that $T$ is transverse to all $C_i$ and it minimizes $\abs{\pi_0(T\cap\cup C_i)}$. If $T$ does not intersect any $C_i$ we have that it is contained in a submanifold homeomorphic to $\Sigma\times I$ which is atoroidal and so $T$ wasn't essential.

		\begin{center}\begin{figure}[h!]
					
						\def\svgwidth{350pt}
						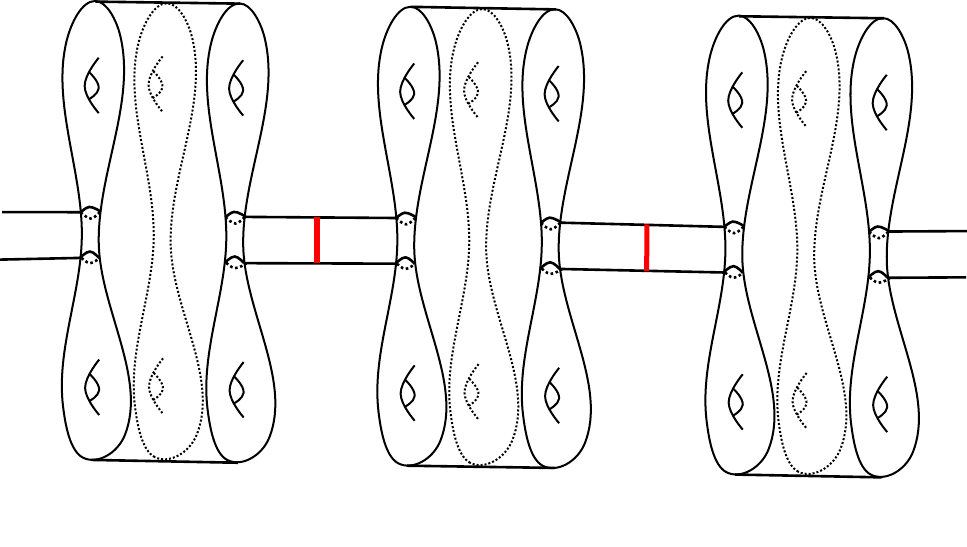
						
						\end{figure}
\end{center}

 Since both $C_i$ and $T$ are incompressible we can isotope $T$ so that the components of the intersection $T\cap C_i$ are essential simple closed curves.\ Thus, $T$ is divided by $\cup_i T\cap C_i$ into finitely many parallel annuli and $T\cap C_i$ are disjoint core curves for $C_i$. Consider $C_k$ such that $T\cap C_k\neq \emp$ and $\forall n\geq k: T\cap C_n=\emp$. Then $T$ cannot intersect $C_k$ in only one component, so it has to come back through $C_k$. Thus, we have an annulus $A\subset T$ that has both boundaries in $C_k$ and is contained in a submanifold of $M_\infty$ homeomorphic to $\Sigma_{k+1}\times I$. The annulus $A$ gives an isotopy between isotopic curves in $\partial \left(\Sigma_{k+1}\times I\right)$ and is therefore boundary parallel. Hence, by an isotopy of $T$ we can reduce $\abs{\pi_0(T\cap\cup C_i)}$ contradicting the fact that it was minimal and non-zero.

We define the submanifold of $M_\infty$ co-bounded by $\Sigma_k$ and $\Sigma_{-k}$ by $M_k$. Since $M_\infty$ is atoroidal so are the $M_k$. Moreover, since the $M_k$ are compact manifolds with infinite $\pi_1$ they are hyperbolizable by Thurston's Hyperbolization Theorem \cite{Kap2001}.

We now want to prove that $M_\infty$ is locally hyperbolic. To do so it suffices to show that given any finitely generated $H\subgroup \pi_1(M_\infty)$ the cover $M_\infty(H)$ corresponding to $H$ factors through a cover $N\twoheadrightarrow M_\infty$ that is hyperbolizable. Let $\gamma_1,\dotsc,\gamma_n\subset M_\infty$ be loops generating $H$. Since the $M_k$ exhaust $M_\infty$ we can find some $k\in\N$ such that $\set{\gamma_i}_{i\leq n} 	\subset M_k$, hence the cover corresponding to $H$ factors through the cover induced by $\pi_1(M_k)$. We now want to show that the cover $M_\infty (k)$ of $M_\infty$ corresponding to $\pi_1(M_k)$	 is hyperbolizable.

Since $\pi:M_\infty\twoheadrightarrow M$ is the infinite cyclic cover of $M$ we have that $M_\infty(k)$ is the same as the cover of $M$ corresponding to $\pi_*(\pi_1(M_k))$. The resolution of the Tameness \cite{AG2004,CG2006} and the Geometrization conjecture \cite{Per2003.1,Per2003.2,Per2003.3} imply the Simon's conjecture, that is: covers of compact irreducible 3-manifolds with finitely generated fundamental groups are tame \cite{Ca2010,Si1976}. Therefore, since $M$ is compact by the Simon's Conjecture we have that $M_\infty(k)$ is tame. The submanifold $M_k\hookrightarrow M_\infty$ lifts homeomorphically to $\widetilde M_k\hookrightarrow M_\infty(k)$. By Whitehead's Theorem \cite{Ha2002} the inclusion is a homotopy equivalence, hence $\widetilde M_k$ forms a Scott core for $M_\infty(k)$. Thus, since $\partial \widetilde M_k$ is incompressible and $M_\infty(k)$ is tame we have that $M_\infty(k)\cong \text{int}(M_k)$ and so it is hyperbolizable.		 
					 
					 \epf

In the infinite cyclic cover $M_\infty$ the essential torus $T$ lifts to a $\pi_1$-injective annulus $A$ that is properly embedded: $A=\gamma\times \R\hookrightarrow M_\infty$ for $\gamma$ the lift of the curve $\alpha\hookrightarrow\Sigma\subset M$.
	 				\brem\label{nosubsurf}
	 				Consider two distinct lifts $\Sigma_i,\Sigma_j$ of the embedded surface $\Sigma\hookrightarrow M$. Then we have that the only essential subsurface of $\Sigma_i$ homotopic to a subsurface of $\Sigma_j$ is a neighbourhood of $\gamma$. This is because by construction the only curve of $\Sigma_i$ homotopic into $\Sigma_j$ is $\gamma$.
					\erem

					 						\bprop\label{esenothyp} The manifold $M_\infty$ is not hyperbolic.
	 						\eprop
	 						\bpf
	 						The manifold $M_\infty$ has two non tame ends $E^\pm$ and the connected components of the complement of a region co-bounded by distinct lifts of $\Sigma$ give neighbourhoods of these ends. Let $A$ be the annulus obtained by the lift of the essential torus $T\hookrightarrow M$. The ends $E^\pm$ of $M_\infty$ are in bijection with the ends $A^\pm$ of the annulus $A$. Let $\gamma$ be a simple closed curve generating $\pi_1(A)$. Denote by $\set{\Sigma_i}_{i\in\Z}\subset M_\infty$ the lifts of $\Sigma\subset M$ and let $\set{\Sigma_i^\pm}_{i\in\Z}$ be the lifts of the punctured tori that form the complement of $\alpha$ in $\Sigma\subset M$. The proof is by contradiction and it will follow by showing that $\gamma$ is neither homotopic to a geodesic in $M_\infty$ nor out a cusp.
\paragraph*{Step 1}  We want to show that the curve $\gamma$ cannot be represented by a hyperbolic element. 
								
								By contradiction assume that $\gamma$ is represented by a hyperbolic element and let $\overline \gamma$ be the unique geodesic representative of $\gamma$ in $M_\infty$. Consider the incompressible embeddings $f_i:\Sigma_2\hookrightarrow  M_\infty$ with $f_i(\Sigma_2)=\Sigma_i$ and let $\gamma_i\subset \Sigma_i$ be the simple closed curve homotopic to $\gamma$.  By picking a 1-vertex triangulation of $\Sigma_i$ where $\gamma_i$ is represented by a preferred edge we can realise each $(f_i,\Sigma_i)$ by a useful-simplicial hyperbolic surface $g_i:S_i\rar M_\infty$ with $g_i(S_i)\simeq \Sigma_i$ (see \cite{Ca1996,Bo1986}). By an abuse of notation we will also use $S_i$ to denote $g_i(S_i)$. Since all the $S_i$ realise $\overline\gamma$ as a geodesic we see the following configuration in $M_\infty$:

	 							\begin{center}\begin{figure}[h!]
	 								\centering
	 								\def\svgwidth{350pt}
	 								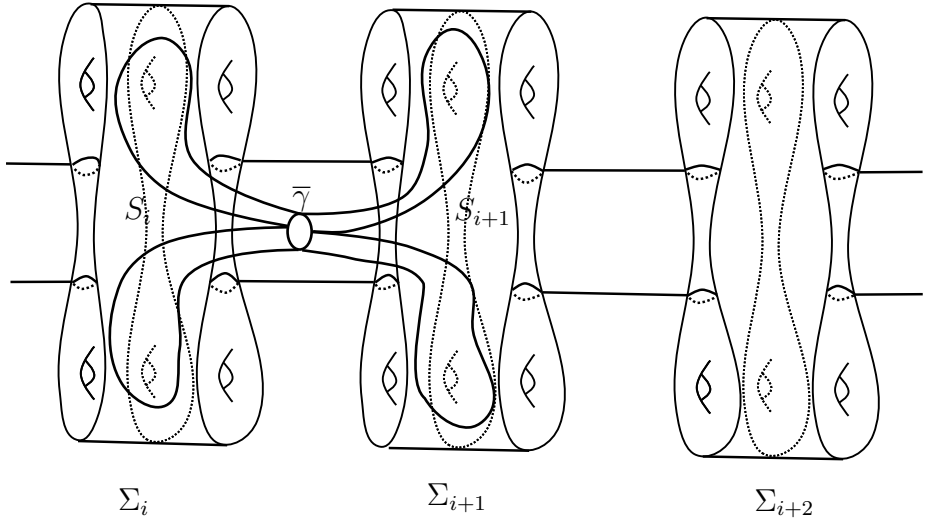
	 									\caption{The simplicial hyperbolic surfaces $S_i$ exiting the ends.}
	 								\end{figure}\end{center}

								On the simplicial hyperbolic surfaces $S_i$ a maximal one-sided collar neighbourhood of $\overline\gamma$ has area bounded by the total area of $S_i$. Since the simplicial hyperbolic surfaces are all genus two by Gauss-Bonnet we have that $ A(S_i)\leq 2 \pi\abs{\chi(S_i)}=4\pi$. Therefore, the radius of a one-sided collar neighbourhood is uniformly bounded by some constant $K=K(\chi(\Sigma_2),\ell(\gamma))<\infty$. Then for $\xi>0$ in the simplicial hyperbolic surface $S_i$ the $K+\xi$ two sided neighbourhood of $\overline\gamma$ is not embedded and contains a 4-punctured sphere. Since simplicial hyperbolic surfaces are $1$-Lipschitz the $4$-punctured sphere is contained in a $K+\xi$ neighbourhood $C$ of $\overline \gamma$, thus it lies in some fixed set $M_h$. Therefore for every $\abs n >h$ we have that $\Sigma_{\pm n}$ has an essential subsurface, homeomorphic to a 4-punctured sphere, homotopic into $\Sigma_{\pm h}$ respectively.\ But this contradicts remark \ref{nosubsurf}.	 

\paragraph*{Step 2} We now show that $\gamma$ cannot be represented by a parabolic element. 

Let $\epsilon>0$ be less then the 3-dimensional Margulis constant $\mu_3$ \cite{BP1992} and let $P$ be a cusp neighbourhood of $\gamma$ such that the horocycle representing $\gamma$ in $\partial P$ has length $\epsilon$. Without loss of generality we can assume that $P$ is contained in the end $E^-$ of $M_\infty$.
	 							 							
Let $\set{\Sigma_i^+}_{i\geq 0}\subset \set{\Sigma_i}_{i\geq 0}$ be the collection of subsurfaces of the $\Sigma_i$ formed by the punctured tori with boundary $\gamma_i$ that are exiting $E^+$. By picking an ideal triangulation of $\Sigma_i$ where the cusp $\gamma_i$ is the only vertex we can realise the embeddings $f_i:\Sigma_i^+\hookrightarrow M_\infty$ by simplicial hyperbolic surfaces $(g_i,S_i^+)$ in which $\gamma_i$ is sent to the cusp \cite{Ca1996,Bo1986}. The $\set{S_i^+}_{i\geq 0}$ are all punctured tori with cusp represented by $\gamma$. 

		\begin{center}\begin{figure}[h!]
	 										\centering
	 										\def\svgwidth{350pt}
	 										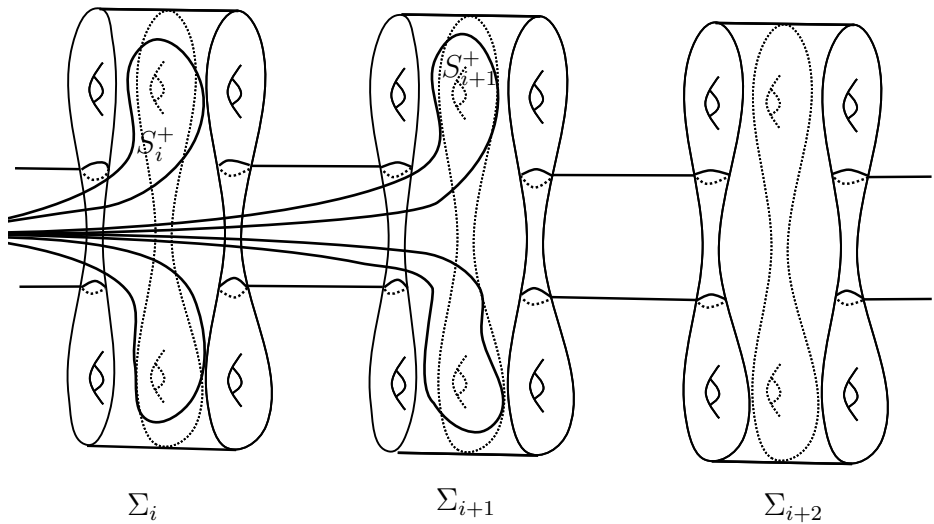
	 										
	 										\caption{The $\epsilon$-thin part is in grey.}
	 									\end{figure}\end{center}

All simplicial hyperbolic surface $S_i^+$ intersects $\partial P$ in a horocycle $f_i(c_i)$ of length $\ell(f_i(c_i))= \epsilon$.\ Therefore, in each $S_i^+$ the horocycle $c_i$ has a a maximally embedded one sided collar whose radius is bounded by some constant $K=K(\epsilon,2\pi)$. Then for $\xi>0$ we have that a $K+\xi$ neighbourhood of $c_i$ in $S_i^+$ has to contain a pair of pants $P_i\subset S_i^+$. Since simplicial hyperbolic surfaces are $1$-Lipschitz the pair of pants of $P_i$ are contained in a $K+\xi$ neighbourhood of $f_i(c_i)$ in $M_\infty$. Thus, the $\Sigma_i$ have pair of pants that are homotopic a uniformly bounded distance from $\partial P$. Let $k\in\N$ be minimal such that $\Sigma_k$ lies outside a $K+\xi$ neighbourhood of $\partial P$. Then for any $j>k$ we have that $\Sigma_j$ has a pair of pants homotopic into $\Sigma_k$ contradicting remark \ref{nosubsurf}. \epf
										
								This concludes the proof of Theorem \ref{agolquest}.		

\nocite{BP1992,Ha2002,He1976,Sh1975,Th1978,Ja1980,MT1998,MT1998,So2006,Bo2006,Ga2006,Si1976}

\thispagestyle{empty}
{\small
\markboth{References}{References}
\bibliographystyle{alpha}
\bibliography{mybib}{}
}

	\bigskip
	\bigskip
	\bigskip
	\bigskip

\noindent Department of Mathematics, Boston College.

\noindent 140 Commonwealth Avenue Chestnut Hill, MA 02467.

\noindent Maloney Hall
\newline \noindent
email: \texttt{cremasch@bc.edu}
	
	\end{document}

%% file: packages_example.tex
\usepackage{latexsym,enumitem}
\usepackage{amssymb}
\usepackage[cp850]{inputenc}
\usepackage{epsfig}
\usepackage{psfrag}
\usepackage{amsthm}
\usepackage{amscd}
\usepackage{amsmath}
\usepackage{amsfonts}
\usepackage{graphics,caption}
\usepackage[all]{xy}
\usepackage{etoolbox}
\patchcmd{\quote}{\rightmargin}{\leftmargin 2em \rightmargin}{}{}
\usepackage[english]{babel}
\usepackage{comment}
\usepackage{color}
\usepackage[colorlinks]{hyperref}
\usepackage[nocompress]{cite} 

\newenvironment{quot}
{
	\vspace{-0.2cm}
	\vspace{0.2cm}
}

\theoremstyle{definition}
\newtheorem{d1}{Definition}

\newenvironment{defin}
{
	\begin{quot}
		\begin{d1}
		}
		{\end{d1}
	\end{quot}

}

\theoremstyle{definition}
\newtheorem{r1}[d1]{Remark}

\newenvironment{rem}
{
	\begin{quot}
		\begin{r1}
		}
		{\end{r1}
	\end{quot}
}

\theoremstyle{definition}
\newtheorem{e1}[d1]{Exercise}

\theoremstyle{definition}
\newtheorem{ese1}[d1]{Example}

\newenvironment{ese}
{
	\begin{quot}
		\begin{ese1}
	}
	{	
		\end{ese1}
	\end{quot}
}

\theoremstyle{definition}

\theoremstyle{definition}
\newtheorem{f2}[d1]{Fact}

\theoremstyle{definition}

\theoremstyle{definition}
\newtheorem*{con2}{Question}

\theoremstyle{definition}
\newtheorem{t1}[d1]{Theorem}

\newenvironment{thm}
{
	\begin{quot}
		\begin{t1}}
		{\end{t1}
	\end{quot}
}

\theoremstyle{definition}
\newtheorem*{T1*}{Theorem}

\newenvironment{teor*}
{
	\begin{quot}
		\begin{T1*}}
		{\end{T1*}
	\end{quot}
}

\newenvironment{dimo}
{\begin{proof}[Proof]
	}
	{\end{proof}}

	\theoremstyle{definition}
	\newtheorem{l1}[d1]{Lemma}
	
	\newenvironment{lem}
	{
		\begin{quot}
			\begin{l1}}
			{\end{l1}
		\end{quot}
	}
	\theoremstyle{definition}
	\newtheorem{p1}[d1]{Proposition}
	
	\newenvironment{prop}
	{
		\begin{quot}
			\begin{p1}}
			{\end{p1}
		\end{quot}
	}
	
	\theoremstyle{definition}
	\newtheorem{c1}[d1]{Corollary}
	
	\newenvironment{cor}
	{
		\begin{quot}
			\begin{c1}}
			{\end{c1}
		\end{quot}
	}

\newtheorem{innercustomthm}{Theorem}
\newenvironment{customthm}[1]
  {\renewcommand\theinnercustomthm{#1}\innercustomthm}
  {\endinnercustomthm}

 \newtheorem*{Theorem*}{Theorem}
 \newtheorem*{Proposition*}{Proposition}
 \newtheorem*{Lemma*}{Lemma}

\let\epsilon\varepsilon
\let\subset\subseteq

\newcommand{\be}{\begin{equation*}}
 \newcommand{\ee}{\end{equation*}}
\newcommand{\bpf}{\begin{dimo}} 
\newcommand{\epf}{\end{dimo}}
\newcommand{\bdefi}{\begin{defin}}
\newcommand{\edefi}{\end{defin}}
\newcommand{\bthm}{\begin{thm}}
\newcommand{\ethm}{\end{thm}}
\newcommand{\blem}{\begin{lem}}
\newcommand{\elem}{\end{lem}}
\newcommand{\bcor}{\begin{cor}}
\newcommand{\ecor}{\end{cor}}

\newcommand{\bprop}{\begin{prop}}
\newcommand{\eprop}{\end{prop}}
\newcommand{\bese}{\begin{ese}} 
\newcommand{\eese}{\end{ese}}
\newcommand{\brem}{\begin{rem}}
\newcommand{\erem}{\end{rem}}

\newcommand{\quotient}[2]{\left.\raisebox{.1em}{$#1\!$}\middle/\raisebox{-.1em}{$#2$}\right.}
\newcommand{\set}[1]{\left\{#1\right\}}	
\DeclareMathOperator{\eqdef}{\doteq\,} 
\DeclareMathOperator{\N}{\mathbb N}			
\newcommand{\abs}[1]{\left\lvert#1\right\rvert}						
\DeclareMathOperator{\R}{\mathbb R}			
\DeclareMathOperator{\emp}{\varnothing} 
\DeclareMathOperator{\Z}{\mathbb Z}			

\newcommand{\bH}{\mathbb H^3}
 \newcommand{\hyp}[1]{\quotient{\bH}{#1}}
  
\newcommand{\subgroup}{\leqslant}

\newcommand{\rar}{\rightarrow} 
\newcommand{\imm}{\looparrowright}



%% file: complex.pdf_tex
\begingroup%
  \makeatletter%
  \providecommand\color[2][]{%
    \errmessage{(Inkscape) Color is used for the text in Inkscape, but the package 'color.sty' is not loaded}%
    \renewcommand\color[2][]{}%
  }%
  \providecommand\transparent[1]{%
    \errmessage{(Inkscape) Transparency is used (non-zero) for the text in Inkscape, but the package 'transparent.sty' is not loaded}%
    \renewcommand\transparent[1]{}%
  }%
  \providecommand\rotatebox[2]{#2}%
  \ifx\svgwidth\undefined%
    \setlength{\unitlength}{377.00653316bp}%
    \ifx\svgscale\undefined%
      \relax%
    \else%
      \setlength{\unitlength}{\unitlength * \real{\svgscale}}%
    \fi%
  \else%
    \setlength{\unitlength}{\svgwidth}%
  \fi%
  \global\let\svgwidth\undefined%
  \global\let\svgscale\undefined%
  \makeatother%
  \begin{picture}(1,1.14811461)%
    \put(0,0){\includegraphics[width=\unitlength,page=1]{complex.pdf}}%
    \put(0.29240554,0.59765929){\color[rgb]{1,0,0}\makebox(0,0)[lb]{\smash{$\alpha\times\set 1$}}}%
    \put(0.2955963,0.46100262){\color[rgb]{1,0,0}\makebox(0,0)[lb]{\smash{$\alpha\times\set 0$}}}%
    \put(0,0){\includegraphics[width=\unitlength,page=2]{complex.pdf}}%
    \put(-0.0018202,0.53445647){\color[rgb]{0,0,0}\makebox(0,0)[lb]{\smash{$\Sigma$}}}%
    \put(0,0){\includegraphics[width=\unitlength,page=3]{complex.pdf}}%
  \end{picture}%
\endgroup%

%% file: cover.pdf_tex
\begingroup%
  \makeatletter%
  \providecommand\color[2][]{%
    \errmessage{(Inkscape) Color is used for the text in Inkscape, but the package 'color.sty' is not loaded}%
    \renewcommand\color[2][]{}%
  }%
  \providecommand\transparent[1]{%
    \errmessage{(Inkscape) Transparency is used (non-zero) for the text in Inkscape, but the package 'transparent.sty' is not loaded}%
    \renewcommand\transparent[1]{}%
  }%
  \providecommand\rotatebox[2]{#2}%
  \ifx\svgwidth\undefined%
    \setlength{\unitlength}{472.3706295bp}%
    \ifx\svgscale\undefined%
      \relax%
    \else%
      \setlength{\unitlength}{\unitlength * \real{\svgscale}}%
    \fi%
  \else%
    \setlength{\unitlength}{\svgwidth}%
  \fi%
  \global\let\svgwidth\undefined%
  \global\let\svgscale\undefined%
  \makeatother%
  \begin{picture}(1,0.49735826)%
    \put(0,0){\includegraphics[width=\unitlength,page=1]{cover.pdf}}%
    \put(0.06420048,0.03634741){\color[rgb]{0,0,0}\makebox(0,0)[lb]{\smash{$\Sigma_i$}}}%
    \put(0.231281,0.0401923){\color[rgb]{0,0,0}\makebox(0,0)[lb]{\smash{$\Sigma_{i+1}$}}}%
    \put(0.40841699,0.03589564){\color[rgb]{0,0,0}\makebox(0,0)[lb]{\smash{$\Sigma_{i+2}$}}}%
    \put(0,0){\includegraphics[width=\unitlength,page=2]{cover.pdf}}%
    \put(0.84676511,0.47261075){\color[rgb]{0,0,0}\makebox(0,0)[lb]{\smash{$\Sigma$}}}%
   
  \end{picture}%
\endgroup%

%% file: horizontalannuli.pdf_tex
\begingroup%
  \makeatletter%
  \providecommand\color[2][]{%
    \errmessage{(Inkscape) Color is used for the text in Inkscape, but the package 'color.sty' is not loaded}%
    \renewcommand\color[2][]{}%
  }%
  \providecommand\transparent[1]{%
    \errmessage{(Inkscape) Transparency is used (non-zero) for the text in Inkscape, but the package 'transparent.sty' is not loaded}%
    \renewcommand\transparent[1]{}%
  }%
  \providecommand\rotatebox[2]{#2}%
  \ifx\svgwidth\undefined%
    \setlength{\unitlength}{278.45553661bp}%
    \ifx\svgscale\undefined%
      \relax%
    \else%
      \setlength{\unitlength}{\unitlength * \real{\svgscale}}%
    \fi%
  \else%
    \setlength{\unitlength}{\svgwidth}%
  \fi%
  \global\let\svgwidth\undefined%
  \global\let\svgscale\undefined%
  \makeatother%
  \begin{picture}(1,0.56502896)%
    \put(0,0){\includegraphics[width=\unitlength,page=1]{horizontalannuli.pdf}}%
    \put(0.12825281,0.0161831){\color[rgb]{0,0,0}\makebox(0,0)[lb]{\smash{$\Sigma_i$}}}%
    \put(0.46208604,0.02141184){\color[rgb]{0,0,0}\makebox(0,0)[lb]{\smash{$\Sigma_{i+1}$}}}%
    \put(0.81601046,0.01556874){\color[rgb]{0,0,0}\makebox(0,0)[lb]{\smash{$\Sigma_{i+2}$}}}%
        \put(0.32681074,0.3564738){\color[rgb]{1,0,0}\makebox(0,0)[lb]{\smash{$C_i$}}}%
    \put(0.66920993,0.34853661){\color[rgb]{1,0,0}\makebox(0,0)[lb]{\smash{$C_{i+1}$}}}%
      \end{picture}%
\endgroup%

%% file: exitingpleat.pdf_tex
\begingroup%
  \makeatletter%
  \providecommand\color[2][]{%
    \errmessage{(Inkscape) Color is used for the text in Inkscape, but the package 'color.sty' is not loaded}%
    \renewcommand\color[2][]{}%
  }%
  \providecommand\transparent[1]{%
    \errmessage{(Inkscape) Transparency is used (non-zero) for the text in Inkscape, but the package 'transparent.sty' is not loaded}%
    \renewcommand\transparent[1]{}%
  }%
  \providecommand\rotatebox[2]{#2}%
  \ifx\svgwidth\undefined%
    \setlength{\unitlength}{278.45553661bp}%
    \ifx\svgscale\undefined%
      \relax%
    \else%
      \setlength{\unitlength}{\unitlength * \real{\svgscale}}%
    \fi%
  \else%
    \setlength{\unitlength}{\svgwidth}%
  \fi%
  \global\let\svgwidth\undefined%
  \global\let\svgscale\undefined%
  \makeatother%
  \begin{picture}(1,0.56502896)%
    \put(0,0){\includegraphics[width=\unitlength,page=1]{exitingpleat.pdf}}%
    \put(0.87562327,0.34197751){\color[rgb]{0,0,0}\makebox(0,0)[lb]{\smash{}}}%
    \put(0.12825281,0.01618308){\color[rgb]{0,0,0}\makebox(0,0)[lb]{\smash{$\Sigma_i$}}}%
    \put(0.46208602,0.02141182){\color[rgb]{0,0,0}\makebox(0,0)[lb]{\smash{$\Sigma_{i+1}$}}}%
    \put(0.81601043,0.01556872){\color[rgb]{0,0,0}\makebox(0,0)[lb]{\smash{$\Sigma_{i+2}$}}}%
    \put(0.31691153,0.34577536){\color[rgb]{0,0,0}\makebox(0,0)[lb]{\smash{$\overline\gamma$}}}%
    \put(0.13487047,0.33116252){\color[rgb]{0,0,0}\makebox(0,0)[lb]{\smash{$S_i$}}}%
    \put(0.49169203,0.32962875){\color[rgb]{0,0,0}\makebox(0,0)[lb]{\smash{$S_{i+1}$}}}%
  \end{picture}%
\endgroup%

%% file: parabolic.pdf_tex
\begingroup%
  \makeatletter%
  \providecommand\color[2][]{%
    \errmessage{(Inkscape) Color is used for the text in Inkscape, but the package 'color.sty' is not loaded}%
    \renewcommand\color[2][]{}%
  }%
  \providecommand\transparent[1]{%
    \errmessage{(Inkscape) Transparency is used (non-zero) for the text in Inkscape, but the package 'transparent.sty' is not loaded}%
    \renewcommand\transparent[1]{}%
  }%
  \providecommand\rotatebox[2]{#2}%
  \ifx\svgwidth\undefined%
    \setlength{\unitlength}{278.45553661bp}%
    \ifx\svgscale\undefined%
      \relax%
    \else%
      \setlength{\unitlength}{\unitlength * \real{\svgscale}}%
    \fi%
  \else%
    \setlength{\unitlength}{\svgwidth}%
  \fi%
  \global\let\svgwidth\undefined%
  \global\let\svgscale\undefined%
  \makeatother%
  \begin{picture}(1,0.56502896)%
    \put(0,0){\includegraphics[width=\unitlength,page=1]{parabolic.pdf}}%
    \put(0.12825281,0.01618308){\color[rgb]{0,0,0}\makebox(0,0)[lb]{\smash{$\Sigma_i$}}}%
    \put(0.46208602,0.02141182){\color[rgb]{0,0,0}\makebox(0,0)[lb]{\smash{$\Sigma_{i+1}$}}}%
    \put(0.81601043,0.01556872){\color[rgb]{0,0,0}\makebox(0,0)[lb]{\smash{$\Sigma_{i+2}$}}}%
     \put(0.87562327,0.34197751){\color[rgb]{0,0,0}\makebox(0,0)[lb]{\smash{}}}%
    \put(0.13791545,0.40779974){\color[rgb]{0,0,0}\makebox(0,0)[lb]{\smash{$S_i^+$}}}%
    \put(0.46822726,0.49013304){\color[rgb]{0,0,0}\makebox(0,0)[lb]{\smash{$S_{i+1}^+$}}}%
     \end{picture}%
\endgroup%